\documentclass[10pt,noamsfonts]{amsart}

\usepackage{verbatim,amsmath,latexsym,amssymb,diagrams,amsfonts}

\theoremstyle{plain}
\newtheorem{thm}{Theorem}[section]

\newtheorem{lem}[thm]{Lemma}

\theoremstyle{definition}
\newtheorem{alg}[thm]{Algorithm}

\newtheorem{ex}[thm]{Example}
\newtheorem{ex-notn}[thm]{Example/Notation}

\newtheorem{fact}[thm]{Fact}

\newtheorem{rem}[thm]{Remark}
\newtheorem{note}[thm]{Note}
\newtheorem{notn}[thm]{Notation}

\def\cone{\operatorname{cone}}	
\def\dR{\operatorname{dR}}

\def\gr{\operatorname{gr}}

\def\im{\operatorname{im}}

\def\proj{\operatorname{Proj}}

\def\spec{\operatorname{Spec}}

\def\var{\operatorname{Var}}

\def\m{{\mathfrak m}}
\def\del{\partial}

\def\A{{\mathbb A}}	
\def\C{{\mathbb C}}

\def\N{{\mathbb N}}

\def\O{{\mathcal O}}	
\def\P{{\mathbb P}}

\def\action{\bullet}
\def\bar#1{\overline{#1}}

\def\into{\hookrightarrow}
\def\Mtwo{{{\em Macaulay} 2\expandafter}}

\def\mylabel#1{\label{#1}}
\def\ignore#1{}

\numberwithin{equation}{section}
\begin{document}

\title[Computing cohomology of complex varieties]  
       {Algorithmic Determination of the Rational Cohomology of
       Complex Varieties via Differential Forms}

\author{Uli Walther}
\address{Department of Mathematics, Purdue University, West 
Lafayette, Indiana 47907}
\email{walther@math.purdue.edu}


\begin{abstract}
We give algorithms for the computation of the algebraic de
Rham cohomology of open and closed algebraic sets inside projective
space or other smooth complex
toric varieties.
The methods, which are based on Gr\"obner basis computations in rings
of differential operators,  can also be used to compute the
cohomology of intersections of smooth closed and open subsets, and in
certain situations the cup-product structure.

We give some examples which were carried out with the help of \Mtwo.
\end{abstract}

\maketitle
\nocite{M2}
\tableofcontents
\section{Introduction}
\mylabel{sec-intro}
The determination of the cohomology of topological spaces has been,
and continues to be, a question of interest going back to (at least)
Poincar\'e. This is documented by the beautiful work of Hopf, Leray,
Serre, and Milnor, to name just a few.  

The advent of reasonably fast
computers brought with it a 
variety Gr\"obner basis driven algorithms performing a multitude of
computations in algebraic and combinatorial settings. This development
did not bypass singular cohomology. In the landmark paper \cite{O-T2}
techniques are presented that compute the dimensions of $H^i(U;\C)$
where $U$ is the complement of an arbitrary algebraic hypersurface in
$\C^n$. These methods were refined in \cite{W-cdrc} in order to deal with
a general Zariski-open set $U\subseteq \C^n$. By \cite{W-cup} one can also
compute the ring structure of $H^*(U;\C)$.

In this note we extend the algorithms from \cite{O-T2, W-cdrc, W-cup} 
to the computation
of cohomology data on more general types of algebraic
sets. These include
\begin{enumerate}
\item singular rational cohomology groups of open sets in projective space,
\item singular rational cohomology of projective varieties,
\item compactly supported rational cohomology of locally closed
varieties in projective space,
\item singular rational cohomology of open subsets of smooth
projective varieties,
\item the ring structure in situation 1.
\end{enumerate}
It follows an overview to the paper.

We shall first give a short review of some known algorithms.
The basic idea of these algorithms is the Grothendieck-Deligne
isomorphism theorem and work by Hartshorne, which assure 
that on complex algebraic spaces  de Rham cohomology can be
computed in the algebraic category, and 
that the singular theory coincides with the
algebraic de Rham theory. At the end of this section we explain some
of the bottlenecks of the algorithms.

We next
consider the special case of projective space. From there we move on
to general open sets in projective space and then, via Alexander duality, to
projective varieties. For open sets, the cup product structure can be
determined. 

Duality can also be used on other spaces, and
that gives access to cohomology of open subsets of smooth projective
varieties, and, as a corollary, compactly supported cohomology of
locally closed sets in projective space.

All the presented methods apply equally well to subvarieties of smooth
toric varieties (as opposed to varieties embedded in projective
space). 
Some of these ideas are expanded in the final
section. 
\begin{notn}
$K$ will be a computable field of characteristic zero contained in
$\C$. Although we will work over $\C$, we shall assume (without
stating this explicitly every time) that all input
data for our algorithms are defined over $K$. This is to guarantee
that we can manipulate the input and recognize vanishing of
expressions with the Turing machine.

Whenever a group is pronounced to be ``finite dimensional'' we
will mean it to be a finite dimensional $\C$-vector space. 
Cosets of elements in a quotient space we usually denote by a bar:
$\bar a$. We write $R_n$ for the ring of polynomials 
$\C[x_1,\ldots,x_n]$, and $D_n$ for the Weyl algebra $\C\langle
x_1,\del_1,\ldots,x_n,\del_n\rangle$. We use multi-index notation:
$x^\alpha\del^\beta$ will mean the monomial
$x_1^{\alpha_1}\cdot\ldots x_n^{\alpha_n}\cdot\del_1^{\beta_1}\cdot\ldots
\cdot\del_n^{\beta_n}$. Also, $|\alpha|$ denotes
in that context $\alpha_1+\ldots +\alpha_n$. 

If $f_0,\ldots,f_r\in R_n$ and $I\subseteq \{0,\ldots,r\}$ we write
$F_I$ for $\prod_{i\in I}f_i$ and $|I|$ for the cardinality of $I$.

If $\phi:K^\bullet\to C^\bullet$ is a chain map of
two chain complexes of modules
over the ring $S$ we write $K^\bullet \cong_S C^\bullet$ if $\phi$ is
a quasi-isomorphism over $S$.

\end{notn}
\section{Review of the affine case}
\mylabel{sec-review}
The purpose of this section is to review an algorithm that leads to
the determination of the singular cohomology groups with
rational coefficients and their ring structure for the complement of
an affine complex variety.
\subsection{}
Let $Y\subseteq X=\spec(R_n)$ 
be defined by the equations $f_0,\ldots,f_r$ in
$R_n$. Then one has associated to $U=X\setminus Y$ a
{\em reduced 
\v Cech complex} 
\begin{eqnarray}
\label{cech}\\
\check C^\bullet=\check C^\bullet(f_0,\ldots,f_r)=
\left(0\to\underbrace{\bigoplus_{|I|=1}R_n[{F_I}^{-1}]}_{\text{degree 0}}\to
\ldots\to\underbrace{\bigoplus_{|I|=r+1}R_n[{F_I}^{-1}]}_{\text{degree }
r}\to 0\right).\nonumber 
\end{eqnarray}
If $U=X$, then we set $\check C^\bullet$ to be the complex
concentrated in degree zero whose entry $\check C^0$ is $R_n$.

One can think of $\check C^\bullet$  as the
appropriate object for various purposes 
that replaces the ring of global sections on $U$ if $U$ is not affine.

$\check C^\bullet(f_0,\ldots,f_r)$ is a complex of (left) $D_n$-modules and 
the maps in the complex are $D_n$-linear
\cite{K2,L-Dmod,W-lcD}.  It makes
therefore sense to speak of bounded complexes $A^\bullet$ of free
finitely generated (left) $D_n$-modules that are
$D_n$-quasi-isomorphic to the \v Cech complex.
\subsection{}
For our purposes we will need a special type of resolutions, 
those that are
{\em $\tilde V_n$-strict} \cite{O-T2,O-T1,W-cdrc}. 
This means that the filtration $\tilde F^\bullet(A^l)$ 
induced by the  grading on $D_n$
defined by $x_i\to 1, \del_i\to -1$ for all $i$ 
is preserved by the maps in the
complex, and that the formation of associated graded objects commutes
with taking homology in $A^\bullet$:
\[
\gr(H^i(A^\bullet[\m_\bullet])\cong H^i(\gr(A^\bullet[\m_\bullet]).
\] 
It is worth pointing out that
this can only be achieved by shifting some of the modules in
$A^\bullet$ appropriately, as in the case of graded resolutions over a
commutative graded ring. 
Even with the shifts, $(A^\bullet[\m_\bullet])$ may not be
graded as the \v Cech complex may not be homogeneous.

It has been shown that for given $f_0,\ldots,f_r$ such a $\tilde V_n$-strict
resolution of $\check C^\bullet=\check C^\bullet(f_0,\ldots,f_r)$ 
is in fact computable. This relies on Gr\"obner basis
techniques, \cite{O-T1, SST, W-cdrc, W-cup}. 

Complexes $A^\bullet[\m_\bullet]\cong_{D_n}\check C^\bullet$ 
that are $\tilde V_n$-strict enjoy a rather
stunning property which we describe now. Consider the Euler operator
$E=x_1\del_1+\ldots +x_n\del_n$. Then $E$ is $\tilde V_n$-homogeneous of
degree 0 and hence acts on $\tilde F^0(A^i[\m_i])/\tilde F^{-1}(A^i[\m_i])$. 
Since
the maps in $A^\bullet[\m_\bullet]$ preserve the
filtration, $E$ acts in fact on $\tilde F^0(H^i)/\tilde F^{-1}(H^i)$ where
$H^i=H^i(A^\bullet)$ with the filtration inherited from $A^i[\m_i]$.
The operator $(-E-n)$ has a minimal nonzero polynomial 
$\tilde b_i(s)$ on this
quotient. (We remark that this holds
not only if $A^\bullet\cong_{D_n}\check C^\bullet$ 
but
more generally  
whenever $A^\bullet$ has {\em holonomic} cohomology, see  \cite{W-cdrc}.)
We write $\tilde b_{A^\bullet[\m_\bullet]}(s)$ for the least common
multiple 
of all these $\tilde b_i(s)$. The polynomial 
$\tilde b_{A^\bullet[\m_\bullet]}(s)$ is
called the {\em $b$-function for integration of
$A^\bullet[\m_\bullet]$ along
$\del_1,\ldots,\del_n$}.  To describe a certain property of $\tilde
b_{A^\bullet[\m_\bullet]}(s)$ we introduce the right module 
$\Omega=D_n/(\del_1,\ldots,\del_n)\cdot D_n$. The functor
$\Omega\otimes^L_{D_n}(-)$ is called {\em integration}.

We say that a cohomology class
$\bar {\bar 1\otimes a}$ in $H^i(\Omega \otimes
_{D_n}A^\bullet[\m_\bullet])$ {\em lives in the  $k$-th level of
the filtration} if $\bar{\bar 1\otimes a}$ has a
representative in $\bar 1\otimes F^{k}(A^i[\m_i])$ but none in 
$\bar 1\otimes F^{k-1}(A^i[\m_i])$. 
The amazing fact is that the roots of $\tilde b_j(s)$ limit the
possible levels of nonzero cohomology classes in
$\Omega \otimes
_{D_n}A^\bullet[\m_\bullet]$. 
Namely, a nonzero class living in
 level $k$ and in cohomological degree $i$ can only occur if 
 $\tilde b_{j}(k)=0$ for some $j\geq i$.

Note that there are only a finite dimensional vector space of
cohomology classes
$\bar{\bar 1\otimes a}$ that live in the $k$-th level of the filtration
because all monomials $x^\alpha\del^\beta$ 
of $\tilde V_n$-degree at most $m-1$
in $D_n[m]$ are right multiples of some  $\del_i$.
From this one can compute the cohomology of
$\Omega\otimes 
_{D_n}A^\bullet[\m_\bullet]$ explicitly because one may simply check
all classes of $\tilde V_n$-degree at most equal to the largest root of 
$\tilde b_{A^\bullet[\m_\bullet]}(s)$.

In a nutshell, this gives the following main steps in an algorithm to
compute the cohomology of $\Omega\otimes
_{D_n}A^\bullet[\m_\bullet]$ (\cite{O-T2,O-T1,W-cdrc}):

\begin{alg}[Integration of the \v Cech complex]~
\mylabel{alg-affine-case}

\noindent {\sc Input}: $f_0,\ldots,f_r\in R_n$, $i\in \N$.

\noindent {\sc Output}: $\dim_\C(H^i(\Omega\otimes 
_{D_n}A^\bullet[\m_\bullet]))$ where $A^\bullet\cong_{D_n}
\check C^\bullet(f_0,\ldots,f_r)$ and $A^\bullet[\m_\bullet]$ is
$D_n$-free and 
$\tilde V_n$-strict. 
\begin{enumerate}
\item Compute a $\tilde V_n$-strict complex $A^\bullet[\m_\bullet]\cong_{D_n}
\check C(f_0,\ldots,f_r)$ (\cite{O-T1,SST,W-cdrc}).
\item Replace each copy of $D_n$ in $A^\bullet$ by 
$\Omega\cong \C[x_1,\ldots,x_n]$. 
\item Find the $b$-functions $\tilde b_i(s)$ for the integration of
$H^i(A^\bullet[\m_\bullet])$ along $\del_1,\ldots,\del_n$, 
and let $k_1$ be the largest
integral root of their product (\cite{O-T1}).
\item Truncate $\Omega\otimes_{D_n} A^\bullet[\m_\bullet]$
to the complex of
finite dimensional $\C$-vector spaces
$\tilde F^{k_1}(\Omega\otimes_{D_n} A^\bullet[\m_\bullet])$  with
$\C$-linear maps.
\item Take the $i$-th cohomology and return its dimension.
\end{enumerate}
End.
\end{alg}
\subsection{}
Now we explain what such a computation has to do with cohomology of
varieties. 
Let $\Omega^\bullet$ be the Koszul complex on $D_n$ induced
by left multiplication by $\del_1,\ldots,\del_n$. Then 
$\Omega\otimes
_{D_n}A^\bullet[\m_\bullet]$ and 
$\Omega^\bullet\otimes_{D_n} \check C^\bullet(f_0,\ldots,f_r)$ are naturally
quasi-isomorphic up to a cohomological shift by $n$. 
This is because $\Omega^\bullet$ 
is a right $D_n$-resolution of $\Omega$ and computing the Tor-functor
can be done by resolving either factor. 

Inspection of the tensor product shows now that it computes algebraic
de Rham cohomology of $U$. The start of this ``inspection'' is
the identification of the complex $\Omega^\bullet\otimes_{D_n}\check
C^\bullet$  in the case $r=0$ 
with the algebraic de Rham functor of \cite{DRCAV} on
$X$ applied to the $\O_X$-module $i_*\O_U=\O_X[f_0^{-1}]$, where
$i:U\into X$. This is why we call $\Omega^\bullet
\otimes_{D_n}\check C^\bullet$ the
{\em algebraic \v Cech-de Rham complex of $U$}.
The
Grothendieck-Deligne comparison theorem and various other ones 
imply that 
$H^i(\Omega^\bullet\otimes_{D_n}\check C^\bullet)\cong_\C
H^{i-n}_{\dR}(U;\C)\cong_\C 
H^{i-n}_{{\rm Sing}}(U;\C)$, the two latter spaces denoting de Rham and
singular cohomology with complex coefficients respectively. 

The essence of the above can be summarized in the following
theorem.
\begin{thm}[de Rham cohomology in affine space \cite{O-T2,W-cdrc}]
\mylabel{thm-affine-dR}
If $f_0,\ldots,f_r$ are given polynomials in $R_n$, then there exists an
algorithm that produces a finite set of cocycles of 
differential forms $\{\omega_{i,i'}\}_{i,i'}$ in
the algebraic \v Cech-de Rham complex $\Omega^\bullet\otimes_{D_n} \check
C^\bullet$ on $U=\C^n\setminus \var(f_0,\ldots,f_r)$ such
that $\{\omega_{i,i'}\}_{i'}$ span $H^i_{\dR}(U;\C)$ for all $i$.
\end{thm}
\begin{proof}
Use Algorithm \ref{alg-affine-case} to obtain a set of generators
(over $\C$) for
the cohomology of $\Omega\otimes_{D_n} A^\bullet[\m_\bullet]$. Then use
Theorem 2.5 of \cite{W-cup} to convert these generators into
cohomology generators for 
$\Omega^\bullet\otimes_{D_n}
\check C^\bullet$
whose elements are  identified with the cochains in the 
algebraic \v Cech-de Rham complex on $U$.
\end{proof}

Now let us give some bibliographical references. 
In \cite{O-T1,W-lcD} algorithms for the presentation of
localizations and more generally the \v Cech complex are discussed. 
In
\cite{O-T2,O-T1,SST,W-cdrc} the $\tilde V_n$-filtration is discussed in varying
detail. 
There it is also explained
how to construct the
complex $A^\bullet[\m_\bullet]$ from $\check C^\bullet$. 
The article \cite{O-T1} gives details
to the computation of the $b$-function and finally \cite{W-cup} shows how
one translates cohomology classes from $\Omega\otimes_{D_n} A^\bullet$ to
classes in $\Omega^\bullet\otimes_{D_n} \check C^\bullet$, thus creating
actual (algebraic) differential forms. 
\subsection{}
It is useful to make some comments about the computational complexity
of the constructions that take place in the execution of Algorithm
\ref{alg-affine-case}. The first major computation is to find a
presentation of the \v Cech complex as a complex of free $D_n$-modules
(Step 1). This computation relies on an algorithm by T.\ Oaku for
determining the Bernstein-Sato polynomial of $f_0\cdot\ldots\cdot
f_r$. Computing this polynomial is quite expensive 
if $\deg(f_0\cdot\ldots\cdot
f_r)>5$. Faster computers will not be of substantial help here because
the complexity of Gr\"obner basis computations usually grows
considerably faster than linearly in the input (the worst possible is
doubly exponentially). Thus, in order to make substantial computational
progress, better algorithms for the Bernstein-Sato polynomial are
needed. 

The next potentially hard step in Algorithm \ref{alg-affine-case} is
to make the 
complex $\tilde V_n$-strict. The author does at this moment 
not know how big a problem this
is. 

A true bottleneck however is the computation of the $b$-functions
$\tilde b_i(s)$, which appears to be somewhat more complex than the
Bernstein-Sato polynomial. But at this time we cannot really make any 
asymptotic
statements. On the positive side, due to the similarity in nature of $\tilde
b_i(s)$ and Bernstein-Sato polynomials one can hope that progress on
one results in progress on the other. 

Step 4 consists of (huge problems in) linear algebra. The author
thinks that this is the least troublesome part of the algorithm, but
whether this is so will much depend on the construction of small $\tilde
V_n$-strict resolutions. 

The algorithms to be described in the sequel use Algorithm 
\ref{alg-affine-case} as a basic building block. None of them involves
computations that make a combinatorial explosion likely to
occur. Unfortunately, however, the current limitations on what
examples can be done with Algorithm \ref{alg-affine-case} restrict us
to rather small examples to illustrate our algorithms. 

\section{Chern classes in projective space}
\mylabel{sec-chern}
In this section we investigate how the cohomology of projective space
can be captured by our formalism. 
It will turn out that it is important to achieve the following.
\begin{lem}
Let  $f_0,\ldots,f_r\in R_n$ be given polynomials, and let
$\{\omega_{i,i'}\}_{i,i'}$ be given cochains of 
algebraic differential forms of
degree $i$ on
$U=X\setminus \var(f_0,\ldots,f_r)$ (i.e.,
$\omega_{i,i'}\in(\Omega^\bullet\otimes_{D_n}\check C^\bullet)^i$). There exists an
algorithm that produces a finite dimensional subcomplex $C^\bullet$ 
of the
algebraic \v Cech-de 
Rham complex $\Omega^\bullet
\otimes_{D_n}\check C^\bullet$ on $U$
such that $C^\bullet\cong_{\C}\Omega^\bullet\otimes_{D_n}\check
C^\bullet$ and $\omega_{i,i'}\in C^i$ $\forall i,i'$.
\end{lem}
\begin{proof}
Let us sketch a proof of the lemma. By Theorem \ref{thm-affine-dR}
it
is possible to find a finite dimensional subcomplex of the algebraic de Rham
complex that captures all the cohomology (namely, just take all the
cohomology generators, with zero differential). However, this may not
include the given forms $\omega_{i,i'}$. Thus, as a first
approximation $C^\bullet_1$ 
of the desired complex $C^\bullet$ we take the union of the
cohomology generators, the given forms $\omega_{i,i'}$ 
and their boundaries $d(\omega_{i,i'})$. 

This is a complex, but the cohomology may be too big. (It is at
least as big as the actual de Rham cohomology but we may have added
extra kernel elements.) We must find forms that reduce the
cohomology. 

By an {\em exhaustion} of a $D_n$-module $M$ we mean a sequence of
$\C$-subspaces $\{D^k(M)\}_{k\in \N}$ such that $D^k(M)\subseteq
D^{k+1}(M)$, $\bigcup_kD^k(M)=M$
and each $D^k(M)$ is finite-dimensional as a $\C$-vector space.
If $M$ is finitely generated over $D_n$ then one may produce an
exhaustion for $M$ from one for $D_n$.
\begin{alg}
\mylabel{alg-lemma-proof}~

\noindent{\sc Input:} $C_1^\bullet$, the subcomplex of
$\Omega^\bullet \otimes_{D_n} \check C^\bullet$ spanned by 
\begin{itemize}
\item the output of the
algorithm of Theorem \ref{thm-affine-dR},
\item forms
$\{\omega_{i,i'}\}_{i,i'}$ with $\omega_{i,i'}\in (\Omega^\bullet
\otimes_{D_n}\check 
C^\bullet)^i$, and their boundaries $\{d(\omega_{i,i'})\}_{i,i'}$.
\end{itemize}
\noindent{\sc Output:} A finite dimensional complex
$C^\bullet\cong_\C\Omega^\bullet\otimes_{D_n}\check C^\bullet$ containing
$C^\bullet_1$ and all $\omega_{i,i'}$.
\begin{enumerate}
\item Initialization: set $l=1$. 
\item Let $i_0=\max\{i:\dim(H^i_{\dR}(U;\C))\not =
\dim(H^i(C_l^\bullet))\}$. If $i_0\le-1$, return $C^\bullet_l$ and
exit. 
\item Let $D^k=D^k(D_n)$ be an exhaustion of $D_n$. For example, let
$D^k=\{x^\alpha\del^\beta:|\alpha+\beta|\le k\}$. 
\item Derive an exhaustion $D^k((\Omega^\bullet \otimes_{D_n} \check
C^\bullet)^{i_0-1})$ of the finitely generated left $D_n$-module
$(\Omega^\bullet \otimes_{D_n} \check
C^\bullet)^{i_0-1}$.
\item For $k=0,1,2,\ldots$ 
test by trial and error whether there is an element in 
$D^k((\Omega^\bullet \otimes_{D_n} \check
C^\bullet)^{i_0-1})$ that maps onto a nonzero element in
$H^{i_0}(C_l^\bullet)$. As soon as such an element is found,  
add it to $C_l^{i_0-1}$, call the enlarged complex $C^\bullet_{l+1}$,
replace $l$ by $l+1$ and move to Step 6. 
\item If $\dim(H^{i_0}_{\dR}(U;\C))\not =
\dim(H^{i_0}(C_l^\bullet))$, reenter at Step 5.
\item Reenter at Step 2.
\end{enumerate}
End.
\end{alg}
\end{proof}

\begin{ex}
Let $r=0$, $n=1$, $x_1=x$  and $f_0=x$. 
We compute $\check C^\bullet=(R_1[x^{-1}])$ positioned in cohomological
degree 0. Moreover, the de Rham cohomology of $U=\C^1\setminus
\var(x)$ is generated by $1$ in degree 0, and $\frac{dx}{x}$ in degree
1.

Suppose  we have the cochain
$\omega_{1,1}=x^{-3}\, dx$ which for some reason we would like to be part
of our finite dimensional subcomplex $C^\bullet$ of
$\Omega^\bullet\otimes_{D_1} \check 
C^\bullet$.

Then $C^\bullet_1$ looks like this:
\[
0\to \C\cdot 1\to \C\cdot\frac{dx}{x}\oplus \C\cdot \frac{dx}{x^3}\to 0,
\]
because $d(x^{-3}dx)=0$.
Clearly
$H^1(C^\bullet_1)$ and $H^1_{\dR}(U;\C)$ do not agree. Thus $i_0=1$ and
we have
to build an exhaustion for $(\Omega^\bullet\otimes_{D_n} \check
C^\bullet)^0=\Omega^0\otimes_{D_n} \check C^0\cong
R_1[x^{-1}]=D_1\action\frac{1}{x}$. 
Take the exhaustion $0\subset D^0\subset D^1\subset\cdots$ 
on $D_1$ given by
$D^k=\C\cdot\{x^a\del^b:a+b\le k\}$. Then $D^0$ is spanned by $\{1\}$,
$D^1$ by $\{1, x ,\del\}$ and $D^2$ by
$\{\del^2,\del,x\del,1,x,x^2\}$. So the exhaustion
$D^\bullet(\Omega^\bullet\otimes_{D_n} \check C^\bullet)^0$ on 
$(\Omega^\bullet\otimes_{D_n} \check C^\bullet)^0
=D_1\action \frac{1}{x}$ in level 0 is spanned by $\{\frac{1}{x}\}$, in
level 1 by $\{\frac{1}{x^2},\frac{1}{x},1\}$ etc.
One
easily sees that the complex $C^\bullet_1
\cup D^0(\Omega^0\otimes_{D_n}\check C^0)\cup d(D^0(\Omega^0\otimes_{D_n}\check
C^0))$ 
has the same first
cohomology as $C_1^\bullet$ while $H^1(C_1^\bullet\cup
D^1(\Omega^0\otimes_{D_n}\check C^0)\cup d(D^1(\Omega^0\otimes_{D_n}\check C^0)))$ 
is one
dimensional. The cause for the drop is 
$\frac{1}{x^2}=\del\action\frac{1}{x}$ with
$d(\frac{1}{x^2})=\frac{-2}{x^3}dx$. 

Set $C^\bullet_2=C^\bullet_1\cup\{\C\cdot\frac{1}{x^2}\}$. 
Then $C^\bullet_2\cong_\C\Omega\otimes_{D_n} \check C^\bullet$
and the algorithm stops.
\end{ex}
\begin{note}
We need to decide linear dependence of a set of given cochains of
differential forms  
in order to make Algorithm \ref{alg-lemma-proof} 
run. This can be
achieved by clearing denominators for example.
\end{note}
For the remainder of the section we shall consider the case of
projective space, which at the same time can be viewed as a warm-up
for general open subsets of $\P^n$,
and as a necessary step for the computation of the cohomology of
(closed) projective varieties.

The major difficulty that arises when  going from affine 
to projective space is
that open sets in projective space are usually not open subsets of an
affine space. This means there is no uniform Weyl algebra
 our computations would be done over. Thus we need to do patching work
and use the Mayer-Vietoris principle.

Recall that the cone $\cone(\phi)$ of a chain map $\phi$ is essentially the
total complex induced by the chain map (\cite{Weibel}, page 18).
\begin{lem}
\mylabel{lem-subquism}
Let $C_1^\bullet$ and $C_2^\bullet$ be two chain complexes of
$S$-modules and let
$K_1^\bullet$ and $K_2^\bullet$ be two subcomplexes which are
quasi-isomorphic:
\[
K_1^\bullet\stackrel{\cong_S}{\into} C_1^\bullet, \quad 
K_2^\bullet\stackrel{\cong_S}{\into} C_2^\bullet. 
\]
Let $\phi:C_1^\bullet\to C^\bullet_2 $ be a chain map that sends
$K_1^\bullet$ into $K^\bullet_2$. Then the cone over $\phi$ is
independent (modulo quasi-isomorphy) of the choice of the complexes:
\[
\cone(K_1^\bullet\stackrel{\phi}{\longrightarrow}K_2^\bullet)
\cong_S
\cone(C_1^\bullet\stackrel{\phi}{\longrightarrow}C_2^\bullet).
\]
\end{lem}
\begin{proof}
The inclusions $K_1^\bullet\stackrel{\cong_S}{\into} C_1^\bullet,  
K_2^\bullet\stackrel{\cong_S}{\into} C_2^\bullet$ induce a map 
$\cone(K_1^\bullet\stackrel{\phi}{\longrightarrow}K_2^\bullet)
\to
\cone(C_1^\bullet\stackrel{\phi}{\longrightarrow}C_2^\bullet)$. To see
that this is a quasi-isomorphism consider the induced map between the long
exact sequences 
\begin{diagram}
H^i(C_2^\bullet)&\rTo&H^i(\cone(C_1^\bullet\stackrel{\phi}{\longrightarrow}C_2^\bullet))&\rTo
&H^i(C_1^\bullet[-1])&\rTo&H^{i+1}(C^\bullet_2)\\
\uTo^{\cong_S}&&\uTo&&\uTo^{\cong_S}&&\uTo^{\cong_S}\\
H^i(K_2^\bullet)&\rTo&H^i(\cone(K_1^\bullet\stackrel{\phi}{\longrightarrow}K_2^\bullet))&\rTo
&H^i(K_1^\bullet[-1])&\rTo&H^{i+1}(K^\bullet_2)
\end{diagram}
and recall the five-lemma.
\end{proof}
Why do we need this lemma? Let us look at $\P^1$. We cover it by the
two open sets $U_1=\spec \C[x]$ and $U_2=\spec \C[x^{-1}]$ which
intersect in $U_{1,2}=\spec \C[x,x^{-1}]$. By Theorem
\ref{thm-affine-dR} 
we know
how to compute differential forms on each of the three open sets that span
the cohomology of the corresponding open set. These would be $\{1_{U_1}\},
\{1_{U_2}\}$, $\{1_{U_{1,2}}$ and 
$(\frac{dx}{x})_{U_{1,2}}\}$. The restriction maps
on the \v Cech-de Rham complex level give us restriction maps $\rho_0^0:
1_{U_1}\to
1_{U_{1,2}}$, $1_{U_2}\to -1_{U_{1,2}}$. (The minus sign is owed to
the general theme of Mayer-Vietoris type complexes.)

One would like to infer that the cohomology of the projective line is
the cohomology of the complex
\small
\begin{diagram}
{\begin{array}{c}\C\cdot 1_{U_1}\\\oplus\\\C\cdot 1_{U_2}\end{array}}\\
       &\rdTo^{\rho^0_0}\\
       &     &\C\cdot 1_{U_{1,2}}&\rTo^0&\C\cdot (\frac{dx}{x})_{U_{1,2}}
\end{diagram}
\normalsize
which equals $\C\cdot(1_{U_1},1_{U_2})$ in degree 0 and
$(\frac{dx}{x})_{U_{1,2}}$ in degree 2.

The reason that this is indeed so is Lemma \ref{lem-subquism}, 
which assures us that
the cohomology of
\small
\begin{diagram}
{\begin{array}{c}
 (\Omega^\bullet_{U_1}\otimes_{D_n} \check C^\bullet_{U_1})^0\\
 \oplus\\
 (\Omega^\bullet_{U_2}\otimes_{D_n} \check C^\bullet_{U_2})^0
\end{array}}
&\rTo&
\begin{array}{c}
 (\Omega^\bullet_{U_1}\otimes_{D_n} \check C^\bullet_{U_1})^1\\
 \oplus\\
 (\Omega^\bullet_{U_2}\otimes_{D_n} \check C^\bullet_{U_2})^1
\end{array}\\
&\rdTo^{\rho^0_0}&&\rdTo^{\rho^0_1}\\
&&(\Omega^\bullet_{U_{1,2}}\otimes_{D_n} \check C^\bullet_{U_{1,2}})^0&\rTo&
(\Omega^\bullet_{U_{1,2}}\otimes_{D_n} \check C^\bullet_{U_{1,2}})^1
\end{diagram}
\normalsize
agrees with the one from the picture above.

\begin{ex-notn}
In this example we investigate projective space
$X=\P^n=\proj(\C[x_0,\ldots,x_n])$.  
Since $X$ is covered by the
$n+1$ open sets $P_j=\spec(\C[x_0,\ldots,x_n,x_j^{-1}]_0)$, the cohomology
and suitable representatives 
can be computed from the combinatorics of this cover.

We write $I$ for a subset of $\{0,\ldots,n\}$ and set $P_I=\bigcap
_{i\in I}P_i$. Then $P_I$ is the set of points in $\P^n$ where
$x_I:=\prod_{i\in I}x_i$ is nonzero. If $i_0=\min_{i\in I}\{i\}$ then 
$P_I\subseteq P_{i_0}\cong \A^n_\C$ is a $(|I|-1)$-fold torus, and its
cohomology is captured by a complex 
\begin{eqnarray}
\label{eqn-dR-open-torus}\\
T^\bullet_I&=&\left(0\to\C\cdot 1\to \bigoplus_{i_0<i\in I} \C\cdot
\frac{d(x_i/x_{i_0})}{x_i/x_{i_0}}\to\cdots\to 
\C\cdot \prod_{i_0<i\in I} \frac{d(x_i/x_{i_0})}{x_i/x_{i_0}}\to 0\right),
\nonumber
\end{eqnarray}
where each differential is zero. Here, the term in cohomological
degree $k$ is
\[
\bigoplus_{{i_0\in J\subseteq I\atop |J|=k+1}}
\C\cdot\prod_{i_0<i\in J}\frac{d(x_i/x_{i_0})}{x_i/x_{i_0}}.
\]
We note that there are several
choices for how to write this complex, because $P_I$ is not only an
open subset of $P_{i_0}$ but also of all other $P_i$ with $i\in
I$. For computing $H^\bullet_{\dR}(\P^n;\C)$ 
we want to glue all these complexes together, for varying $I$. Then we
need to translate differential forms from the chart $x_j\not =0 $ to
those on the chart $x_{j'}\not =0$. 
\begin{fact}
\mylabel{fact-translation}
The conversion of differential forms on $P_j\cap P_{j'}$ 
from the $j$-chart $P_j$
to the $x_{j'}$-chart $P_{j'}$ is obtained as
follows.
\begin{itemize}
\item 
$f(x_0/x_j,\ldots,\widehat{x_j/x_j},\ldots,x_n/x_j)=
\frac{f(x_0/x_{j'},\ldots,\widehat{x_{j'}/x_{j'}},\ldots,x_n/x_{j'})} 
{({x_j}/{x_j'})^{\deg(f)}}$,
\item
for all $i\not = j,j'$ we have
$d(x_i/x_j)=\frac{d(x_i/x_{j'})}{(x_j/x_{j'})}-
\frac{(x_i/x_{j'})d(x_j/x_{j'})}{(x_j/x_{j'})^{2}}$,
\item $d(x_{j'}/x_j)=\frac{-d(x_j/x_{j'})}{(x_j/x_{j'})^{2}}$.
\end{itemize}
All these expressions are regular on  $P_j\cap
P_{j'}=\spec(K[x_0,\ldots,x_n,x_j^{-1},x_{j'}^{-1}]_0)$. 
\end{fact}

Note that by Fact \ref{fact-translation}, the space
$T^k_I$ is formally invariant under
the change of charts $j\to j'$. 

Since the de Rham cohomology of $\P^n$ is the cohomology of the total
complex of the global sections of the \v Cech-de Rham complexes on the open
tori, we conclude by Lemma \ref{lem-subquism} that the cohomology is captured by the
total complex of the complexes (\ref{eqn-dR-open-torus}).

For $n=2$ this looks like this:
\small
\begin{diagram}
{\begin{array}{c}\C\cdot 1\\\oplus\\\C\cdot 1\\\oplus\\\C\cdot 1\end{array}}\\
&\rdTo^{\rho_0^0\cdot}\\
&&{\begin{array}{c}\C\cdot 1\\\oplus\\\C\cdot 1\\\oplus\\\C\cdot 1\end{array}}&\rTo&
  {\begin{array}{c}
  \C\cdot \frac{d(x/y)}{x/y}\\\oplus\\
  \C\cdot \frac{d(x/z)}{x/z}\\\oplus\\
  \C\cdot \frac{d(y/z)}{y/z}\end{array}} &\\
&&&\rdTo_{\rho^0_1\cdot}&&\rdTo^{\rho^1_1\cdot}\\
&&&&\C\cdot 1&\rTo&\C\cdot \frac{d(x/z)}{x/z}\oplus
\C\cdot \frac{d(y/z)}{y/z}&\rTo&\C\cdot \frac{d(x/z)\wedge d(y/z)}{(x/z)(y/z)}
\end{diagram}
\normalsize
Here, the rows correspond to the subcomplexes
(\ref{eqn-dR-open-torus})
of the algebraic \v Cech-de Rham
complexes on $P_1,P_2,P_3$ (top block), $P_{1,2}, P_{1,3}, P_{2,3}$
(middle block) and $P_{1,2,3}$ (bottom line). 
All maps are zero in
horizontal direction, and 
$\rho^0_0=\left(\begin{array}{ccc}1&1&0\\-1&0&1\\0&-1&-1\end{array}\right)$,
$\rho^0_1=(1,-1,1)$ and
$\rho^1_1=\left(\begin{array}{ccc}1&-1&0\\-1&0&1\end{array}\right)$.

One can see (here, as well as in general) that the cohomology of
$\P^n$ is one-dimensional in even degree $2k$, generated by the
$k$-cocycle of $k$-forms 
\begin{eqnarray}
\label{chern-classes}
\label{eqn-std-open-sets}
c_k&=&\sum_{I=\{i_0<\ldots<i_k\}\subseteq\{0,\ldots,n\}} 
\prod_{i_0<i\in
I}\,\,\frac{d(x_i/x_{i_0})}{x_i/x_{i_0}}.   
\end{eqnarray}
where the displayed summand is 
defined on the $k+1$-fold intersection
$P_I=\P^n\setminus\var(x_{i_0}\cdot\ldots x_{i_k})$.

The $c_k$ are, up to a
constant, the Chern classes of projective space.
\end{ex-notn}

The example gives an indication how more general open sets
will be attacked, namely by an open cover, used in conjunction with 
 the translation
formul\ae\ from Fact \ref{fact-translation}.

\section{Open sets in projective space}
\mylabel{sec-general}
Let $U$ be the open set of $\P^n$ defined by the non-vanishing
of $f_0,\ldots, f_r\in S=K[x_0,\ldots,x_n]$. 
In this section we describe an algorithm to compute the rational cohomology of
 $U$ using the algebraic \v Cech-de Rham complex
on an open cover by the sets $U_I=U\cap P_I$ where the $P_I$
are the open sets from \ref{eqn-std-open-sets} covering $\P^n$.

Let $I\subseteq \{0,\ldots,n\}$ and $j=\min_{i\in I}(i)$.
We consider $U_I$ as the open set in
$P_j=\spec(\C[\frac{x_0}{x_j},\ldots,\frac{x_n}{x_j}])$ whose
complement is the variety of 
$\{
\frac{x_I}{{x_j}^{|I|}}f_i(\frac{x_0}{x_j},\ldots,\frac{x_n}{x_j})\}_{i=1}^r$.
We denote by $D_I$ the Weyl algebra associated to the ring
$\C[\frac{x_0}{x_j},\ldots,\frac{x_n}{x_j}]$ and by $\Omega^\bullet_I$
the Koszul complex of right $D_I$-modules induced by
$\del_{x_0/x_j},\ldots,\del_{x_n/x_j}$. 

Algorithm \ref{alg-affine-case} in combination with Theorem
\ref{thm-affine-dR}  
produces
for each such $I$
a finite number of cocycles of
differential forms  which generate the algebraic de Rham cohomology 
$H^\bullet_{\dR}(U_I;\C)$.

We can think of these classes for fixed $I$ as a
subcomplex of the algebraic \v Cech-de Rham
complex on $U_I$ with
zero differential, having the same cohomology as the whole algebraic
de Rham 
complex  on $U_I$. Our goal is to glue these complexes according to
the open cover, and compute cohomology.

Unfortunately, the natural maps of differential
forms induced by the inclusions $U\cap P_I\cap 
P_j\into U\cap P_I$ may not be carried by these
subcomplexes. Since we need this to happen in order to form a 
total complex from the subcomplexes and to use Lemma
\ref{lem-subquism} we need to enlarge the subcomplexes suitably.

The strategy is to start with the complex on $U\cap P_j$, $0\le j\le
n$, and work our
way up to higher and higher intersections. What we need to achieve is
a set of finite dimensional complexes $C_I^\bullet$ on $U\cap
P_I$ such that if $c\in C_I^\bullet$
then its natural image in the algebraic \v Cech-de Rham complex on $U\cap P_I$ 
is in $C^\bullet_{I\cup j}$ for all $I,j$. (This natural
image is of course for each differential form given by exactly the
same form, considered as a form on an open subset.)

Let us give an outline for how to do one such step. Take
$C^\bullet_I$ and $C^\bullet_{I\cup j}$ 
where the former was obtained from the
integration if $|I|=1$ and from the inductive step otherwise, while
the latter 
comes from Theorem \ref{thm-affine-dR}.

Execute  Algorithm
\ref{alg-lemma-proof} with the following input and output variables. 
For $C_1^\bullet$ we take  $C^\bullet_{I\cup j}$.
The set $\{\omega_{i,i'}\}_{i'}$ is for
each $i$ a set of vector space generators for $C_I^i$. The output
$C^\bullet$ is
quasi-isomorphic to $C^\bullet_{I\cup j}$, contains $C^\bullet_I$, 
and replaces the old
(input) complex $C^\bullet_{I\cup j}$. 

Iterating over $|I|$ from 1 to $n$ we get a collection 
of finite dimensional complexes $C^\bullet_I$ of differential
forms whose $k$-th cohomology is exactly $H^k_{\dR}(U_{I};\C)$,
and $C^\bullet_I\into\Omega^\bullet_I\otimes_{D_I} \check C^\bullet_I\to
\Omega^\bullet_{I\cup j}\otimes_{D_{I\cup j}} \check C^\bullet_{I\cup
j}$ factors  
through $C^\bullet_{I\cup j}\into \Omega^\bullet_{I\cup
j}\otimes_{D_{I\cup j}}\check C^\bullet_{I\cup j}$.
By Lemma \ref{lem-subquism} the total complex composed of the
complexes $C^\bullet_I$ is 
quasi-isomorphic to the algebraic 
\v Cech-de Rham complex on $U$ relative to the
cover $U_I$.  We hence have
\begin{alg}[Cohomology of open sets]
\mylabel{alg-proj-open}~

\noindent{\sc Input:} Homogeneous polynomials $f_0,\ldots,f_r$ in
$K[x_0,\ldots,x_n]$.

\noindent{\sc Output:} The cohomology groups of $U=\P^n_\C\setminus
\var(f_0,\dots,f_r)$. 
\begin{enumerate}
\item For each $I\subseteq \{0,\ldots,n\}$ 
compute a finite dimensional 
subcomplex $C^\bullet_I$ of the \v Cech-de Rham complex on
$U_I=U\cap P_I$ with $C^\bullet_I\cong_\C\Omega^\bullet_I\otimes_{D_I}\check
C^\bullet_I$ (Theorem \ref{thm-affine-dR}).
\item For $k=1,2,3,\ldots,n$ do
\begin{itemize}
\item for all $|I|=k, I\subseteq \{0,\ldots,n\}$ do
\begin{itemize}
\item for all $j\in \{0,\ldots,n\}\setminus I$ do

run Algorithm \ref{alg-lemma-proof} with
\begin{itemize}
\item Input:
\begin{itemize}
\item $\{\omega_{i,i'}\}_{i'}:=$ a set of vector space generators for
$C^i_I$;
\item $C^\bullet_1:= C^\bullet_{I\cup j}$
\end{itemize}
\item Output $C^\bullet$ replacing $C^\bullet_{I\cup j}$.
\end{itemize}
\end{itemize}
\end{itemize}
\item Set up the total complex $C^\bullet_U$ induced by the maps
$C^\bullet_I\into C^\bullet_{I\cup j}$.
\item Compute the cohomology of $C^\bullet_U$ which equals 
the singular (or de Rham) cohomology of
$\P^n_\C\setminus \var(f_0,\ldots,f_r)$.
\end{enumerate}
End.
\end{alg}
The elements, and hence the cohomology, of $C^\bullet_U$ are \v Cech
cochains  of differential forms for the cover $U=\bigcup_{0\le j\le
n}U_j$. 

We remark that similarly to the affine case the complex $C_U^\bullet$
carries not quite enough information to compute the cup product
structure of $U$, but that also like in the affine case this can be
fixed by further enlarging $C_U^\bullet$: 
\begin{alg}[Cup products on open sets]
\mylabel{alg-cup-proj-open}~

\noindent{\sc Input:} Homogeneous polynomials $f_0,\ldots,f_r$ in
$K[x_0,\ldots,x_n]$.

\noindent{\sc Output:} A multiplication table for
$H^\bullet_{\dR}(U;\C)$.

\begin{enumerate}
\item Run Algorithm \ref{alg-proj-open} to get the complex
$C^\bullet_U$. 
\item Compute explicit generators (cocycles of differential forms) 
for the cohomology of $C_U^\bullet$.
\item Multiply these forms in the \v Cech-de Rham complex on $U$
according to the usual rules for multiplying \v Cech cochains, see for
example \cite{W-cup}, Theorem 4.1.
\item Enlarge $C^\bullet_U$ so that it contains all these products
(using Algorithm \ref{alg-lemma-proof}).
\item Determine a presentation of the cosets of 
the products in terms of the chosen
representatives for the cohomology of $C^\bullet_U$ to get a
multiplication table.
\end{enumerate}
End.
\end{alg}
\section{An example}
\mylabel{sec-example}
In this section we will go through one example in detail: the curve
$C=\var(x^2+yz)$ in $\P^2$. This is of course a rather specific example, but
more interesting examples are too large to be 
 useful for an illustration of the general technique (and, as outlined
in the introduction, examples of substantial interest are out of reach
at the moment).

On the three coordinate patches of $\P^2$, the complement of $C$ is
given by the non-vanishing of $1+(z/x)(y/x)$, $(x/y)^2+z/y$ and $(x/z)^2+y/z$. 
We shall call $U_1,U_2,U_3$ 
the corresponding coordinate patches of
$\P^2\setminus C$, and $U_{1,2}$, $U_{1,3}$,
$U_{2,3}$ and $U_{1,2,3}$ their intersections.

First we determine a finite set of
(exact) differential forms on each of the $U_I$ such that the inclusion of
the complex $C^\bullet_I$ 
generated by these forms 
(with trivial differential) 
into the algebraic \v Cech-de Rham complex $\Omega^\bullet_I\otimes
_{D_I}\check C^\bullet_I$ is
a quasi-isomorphism. Later we shall consider the natural maps
$\Omega^\bullet_I\otimes_{D_I}\check C^\bullet_I\to 
\Omega^\bullet_{I\cup j}\otimes_{D_{I\cup j}}\check C^\bullet_{I\cup j}$
obtained from the inclusions $U_{I\cup
j}\into U_I$.

With \Mtwo\ one computes that the de Rham cohomology groups of the
various $U_I$ are generated by the following elements. 
\small
\[
\begin{array}{|c|c|ccc|}
\hline
U_1&H^0&&1&\\\cline{2-5}
   &H^1&e_{1:1}&=&
        {\frac{(y/x)(z/x)^2\,d(y/x)+(y/x)^2(z/x)\,d(z/x)}{(1+(y/x)(z/x))^2}}\\
   &H^2&t_{1:1}&=&\frac{(y/x)(z/x)\,d(y/x)\,d(z/x)}{(1+(y/x)(z/x))^2}\\\cline{2-5}
\hline
U_2&H^0&&1&\\\cline{2-5}
   &H^1&e_{2:1}&=&\frac{d(z/y)+2(x/y)\,d(x/y)}{(x/y)^2+(z/y)}\\
\hline
U_3&H^0&&1&\\\cline{2-5}
   &H^1&e_{3:1}&=&\frac{d(y/z)+2(x/z)\,d(x/z)}{(x/z)^2+(y/z)}\\
\hline
U_{1,2}&H^0&&1&\\\cline{2-5}
       &H^1&e_{1,2:1}&=&\frac{(x/y)(z/y)^2\,d(x/y)-
               \frac{1}{2}(x/y)^2(z/y)\,d(z/y)}
            {(((x/y)^2+(z/y))(x/y))^2}\\
       &   &e_{1,2:2}&=&\frac{d(x/y)}{(x/y)}\\\cline{2-5}
       &H^2&t_{1,2:1}&=&\frac{(x/y)(z/y)\,d(x/y)\,d(z/y)}{(((x/y)^2+(z/y))(x/y))^2}\\
\hline
U_{1,3}&H^0&&1&\\\cline{2-5}
       &H^1&e_{1,3:1}&=&\frac{(x/z)(y/z)^2\,d(x/z)-\frac{1}{2}(x/z)^2(y/z)\,d(y/z)}
           {(((x/z)^2+(y/z))(x/z))^2}\\
       &   &e_{1,3:2}&=&\frac{d(x/z)}{(x/z)}\\\cline{2-5}
       &H^2&t_{1,3:1}&=&\frac{(x/z)(y/z)\,d(x/z)\,d(y/z)}
           {(((x/z)^2+(y/z))(x/z))^2}\\
\hline
U_{3,2}&H^0&&1&\\\cline{2-5}
       &H^1&e_{3,2:1}&=&\frac{-2(y/z)^2(x/z)^3\,d(x/z)+(y/z)(x/z)^4\,d(y/z)}
           {(((x/z)^2+(y/z))(y/z))^2}\\
       &   &e_{3,2:2}&=&\frac{d(y/z)}{(y/z)}\\\cline{2-5}
       &H^2&t_{3,2:1}&=&
            \frac{(y/z)(x/z)^3\,d(x/z)\,d(y/z)}{(((x/z)^2+(y/z))(y/z))^2}\\
\hline
U_{1,2,3}&H^0&&1&\\\cline{2-5}
         &H^1&e_{1,2,3:1}&=&\frac{-2(y/z)^2(x/z)^5\,d(x/z)+(y/z)(x/z)^6}
         {(((x/z)^2+(y/z))(x/z)(y/z))^2}\\
         &   &e_{1,2,3:2}&=&\frac{-2(y/z)^4(x/z)\,d(x/z)+(y/z)^3(x/z)^2\,d(y/z)}
         {(((x/z)^2+(y/z))(x/z)(y/z))^2}\\
         &   &e_{1,2,3:3}&=&\frac{d(y/z)}{(y/z)}\\\cline{2-5}
         &H^2&t_{1,2,3:1}&=&\frac{(y/z)(x/z)^5\,d(x/z)\,d(y/z)}
             {(((x/z)^2+(y/z))(x/z)(y/z))^2}\\
         &   &t_{1,2,3:2}&=&\frac{(y/z)^3(x/z)\,d(x/z)\,d(y/z)}
             {(((x/z)^2+(y/z))(x/z)(y/z))^2}\\
\hline
\end{array}
\]
\normalsize
In this table, $e_{I:k}$ is the $k$-th generator of $H^1(U_I;\C)$ while
$t_{I:k}$ is the $k$-th generator of $H^2(U_I;\C)$. 
For example,
the commands for $U_3$ are 
\begin{verbatim}
load "../m2/Dloadfile.m2"
R=QQ[s,t]   -- s=x/z, t=y/z
f=s^2+t
deRhamAll(f)
\end{verbatim}
The first line loads the $D$-module library \cite{M2D}.
From \ref{fact-translation},
\begin{eqnarray*}
d(x/y)&=&(y/z)^{-1}d(x/z)-(x/z)(y/z)^{-2}d(y/z),\\
d(z/y)&=&-(y/z)^{-2}d(y/z).
\end{eqnarray*} 
This holds of course for any permutation of the variables $x,y,z$ as
well.
Set
$g_{y,z}=\frac{(x/z)^2(y/z)}{((x/z)^2+(y/z))(y/z)}$, $g_{x,z}=
\frac{(x/z)(y/z)}{((x/z)^2+(y/z))(x/z)}$,
$g_{x,y}=\frac{(x/y)(z/y)}{((x/y)^2+(z/y))(x/y)}$ and note that
$g_{x,y}=g_{x,z}=1-g_{y,z}$. 

With these rules and abbreviations one computes the following identifications
representing the maps from 0-cochains to 1-cochains and from 1-cochains
to 2-cochains.
\small
\begin{eqnarray*}
e_{3:1}&=&-e_{2,3:1}+e_{2,3:2}+d(g_{y,z})=2e_{1,3:2}-2e_{1,3:1}+d(g_{x,z})\\
e_{2:1}&=&-e_{2,3:1}-e_{2,3:2}+d(g_{y,z})=2e_{1,2:2}-2e_{1,2:1}+d(g_{x,y})\\
e_{1:1}&=&-2e_{1,3:1}=-2e_{1,2:1}\\
t_{1:1}&=&t_{1,3:1}=-t_{1,2:1}\\
e_{2,3:1}&=&e_{1,2,3:1}\\
e_{2,3:2}&=&e_{1,2,3:3}\\
t_{2,3:1}&=&t_{1,2,3:1}\\
e_{1,3:1}&=&-\frac{1}{2}e_{1,2,3:2}\\
e_{1,3:2}&=&\frac{1}{2}e_{1,2,3:1}-\frac{1}{2}e_{1,2,3:2}+
            d(g_{y,z})\\
t_{1,3:1}&=&-t_{1,2,3:2}\\
e_{1,2:1}&=&-2e_{1,2,3:2}\\
e_{1,2:2}&=&-\frac{1}{2}e_{1,2,3:1}-\frac{1}{2}e_{1,2,3:2}+e_{1,2,3:3}
+d(g_{y,z}).
\end{eqnarray*}
\normalsize
Then the following forms generate finite dimensional complexes
$C^\bullet_I$ that are quasi-isomorphic by the inclusion 
to the \v Cech-de Rham complex on
$U_I$.
\small
\[
\begin{array}{|c|ccc|}
\hline
U_1&1&e_{1:1}&t_{1:1}\\
U_2&1&e_{2:1}&\\
U_3&1&e_{3:1}&\\
U_{1,2}&1,
g_{x,y}&e_{1,2:1},e_{1,2:2},
d(g_{x,y})&t_{1,2:1}\\
U_{1,3}&1
,g_{x,z}& e_{1,3:1},e_{1,3:2},
d(g_{x,z})&t_{1,3:1}\\
U_{2,3}&1,
g_{y,z}& e_{2,3:1}, e_{2,3:2},
d(g_{y,z})& t_{2,3:1}\\
U_{1,2,3}&1,g_{y,z}& e_{1,2,3:1},
e_{1,2,3:2}, e_{1,2,3:3}, d(g_{y,z})& 
t_{1,2,3:1}, t_{1,2,3:2}\\\hline
\end{array}
\]
\normalsize
The forms $g_{x,y}, g_{x,z}$ and $g_{y,z}$ and their boundaries are
needed to assure that $C^\bullet_I\subseteq C^\bullet_{I\cup j}$ for
all $I\subseteq\{1,2,3\}$ and all $j\in\{1,2,3\}\setminus I$. 

By Lemma \ref{lem-subquism} the \v Cech-de Rham complex on
$U=\P^2\setminus C$ is quasi-isomorphic to the total complex composed of
the $C_I^\bullet$. 

Let $U_{(1)}=\{U_1,U_2,U_3\}$, $U_{(2)}=\{U_{1,2},
U_{2,3}, U_{1,3}\}$. Then the total complex $C^\bullet_U$ 
made from the
$C^\bullet_I$ has
\begin{itemize}
\item three terms in degree zero (the three constants on $U_{(1)}$),
\item nine terms in degree 1 (three constants from
$U_{(2)}$, three $H^1$-generators from $U_{(1)}$ and
the three 0-forms $g_{x,y}$, $g_{x,z}$, $g_{y,z}$ on $U_{(2)}$),
\item twelve terms in degree 2 (one $H^2$ generator from $U_3$, six
$H^1$-generators from $U_{(2)}$, the constants from
$U_{1,2,3}$, the differentials of the extra 0-forms on $U_{(1)}$ 
and an additional 0-form on $U_{1,2,3}$),
\item seven terms in degree 3 (three $H^2$-generators from $U_{(2)}$, three
$H^1$-generators from $U_{1,2,3}$ and the differential of the
additional 0-form on $U_{1,2,3}$),
\item two terms in degree 4 (the $H^2$-generators on $U_{1,2,3}$).
\end{itemize}
Since this is a finite-dimensional complex, we can
compute its cohomology by linear algebra. This  
determines \v Cech-de Rham cochains of 
differential forms in the \v
Cech-de Rham complex on $U$ that carry the de Rham cohomology of $U$.
With \Mtwo\ again one
computes the cohomology of this complex to be zero in all degrees but
in $H^0$ where the cohomology is isomorphic to $\C$, generated by the
cochain $(1,1,1)$ of 0-forms on $U_{(1)}$.
\section{Closed and locally closed subsets}
\mylabel{sec-compact}
\subsection{Closed varieties in $\P^n$}
\mylabel{subsec-closed}
In this subsection we consider what information can be obtained of the
cohomology of the closed sets $Y=\var(f_0,\ldots,f_r)$.

We first note that there is a long exact sequence of sheaf cohomology 
\[
\cdots\to H^k_Y(\P^n;\C)\to H^k(\P^n;\C)\to H^k(U;\C)\to
H^{k+1}_Y(\P^n;\C)\to\cdots
\]
Here $\C$ denotes the constant sheaf. 
Furthermore, since $\P^n$ is a manifold of dimension $2n$, 
we can use Alexander duality \cite{Iversen}.
Hence
\[
H^k_Y(\P^n,\C)\cong H^{2n-k}_c(Y;\C)^*,
\]
the latter denoting the vector space dual of
cohomology with compact supports.
Since $Y$ is a compact space however, cohomology with compact supports
agrees with the usual cohomology. Considering the structure of the
cohomology groups on $\P^n$ (zero in odd degree) 
there are exact sequences
\begin{eqnarray*}
0\to H^{2k-1}(U;\C)\to H^{2n-2k}(Y;\C)^*\to
H^{2k}(\P^n;\C)\phantom{xxxxxxxxxx}\\
\to
H^{2k}(U;\C)\to H^{2n-2k-1}(Y;\C)^*\to 0
\end{eqnarray*}
for all $i>0$ and an exact sequence 
\[
0\to H^{2n}(Y;\C)^*\to H^0(\P^n;\C)\to H^0(U;\C)\to H^{2n-1}(Y;\C)^*\to 0
\]
where the star denotes the vector space dual.
It is not
 hard to understand the maps $H^k(\P^n;\C)\to H^k(U;\C)$
algorithmically. In Section \ref{sec-chern} we found generators for the
cohomology of $\P^n$. Since the inclusion $U\to \P^n$ induces the maps
$H^k(\P^n;\C)\to H^k(U;\C)$, the forms $c_k\in H^{2k}(\P^n;\C)$ 
 are simply interpreted as forms on $U$. In order to find the kernel
 and the cokernel of $H^k(\P^n;\C)\to H^k(U;\C)$
 it is sufficient to find a subcomplex $C^\bullet\cong_\C C^\bullet_U$ 
of the algebraic \v
 Cech-de Rham complex on $U=\bigcup(P_i\cap U)$ that contains each
 $c_k$, because of Lemma \ref{lem-subquism}. Such a complex can be
constructed from Algorithm \ref{alg-lemma-proof}.

Hence we have the following algorithm:
\begin{alg}[Cohomology of projective varieties]
\mylabel{alg-proj-coh}~

\noindent{\sc Input:} Homogeneous polynomials $f_0,\ldots,f_r$ in
$\C[x_0,\ldots,x_n]$.

\noindent{\sc Output:} The cohomology groups of $
\var(f_0,\dots,f_r)$ in $\P^n_\C$. 

Let $U=\P^n\setminus\var(f_0,\ldots,f_r)$.
\begin{enumerate}
\item Compute $C^\bullet_U$, a finite dimensional 
complex of differential forms on
$U=\bigcup_{i=0}^r (U\cap P_i)$ from Algorithm
\ref{alg-proj-open} that computes the de Rham cohomology of $U$.
\item Use Algorithm \ref{alg-lemma-proof} to enlarge 
$C^\bullet_U$ so that it contains for all $k$ the cocycles $c_k$
from (\ref{chern-classes}).
\item $H^{2n-2k-1}(Y;\C)$ is isomorphic to 
the cokernel of the map 
\[
H^{2k}(\P^n;\C)\to
H^{2k}(U;\C). 
\]
The dimension of this space agrees with
$\dim_\C H^{2k}(U;\C)$ if $c_k$ represents the zero class in
$H^{2k}(U;\C)$; else it is $\dim_\C
H^{2k}(U;\C)-1$. The vanishing of 
$c_k$ in $H^{2k}(U;\C)$ is equivalent to $c_k$
being an image in $C^\bullet_U$.
\item $\dim_\C H^{2n-2k}(Y;\C)$ equals 
$\dim_\C H^{2k-1}(U;\C)$ if $c_k=0$ in $H^{2k}(U;\C)$ and $\dim_\C
H^{2k-1}(U;\C)+1$ otherwise.
\end{enumerate}
End.
\end{alg}
\begin{rem}
Since Chern classes are preserved under
pullbacks, 
if $U$
is some open set in $\P^n$ then the images on $U$ 
of the generators for
$H^{2k}(\P^n;\C)$ from (\ref{eqn-dR-open-torus})
are the Chern classes of $U$. This shows how one can
determine vanishing of the rational Chern classes on $U$.
$Y$ and $U$ have the same cohomology dimensions except
for a difference of 1 or $-1$. This difference is dictated by the
vanishing of the Chern classes of $U$. 
\end{rem}
\begin{ex}
\mylabel{veronese-closed}
Consider the variety defined by $x^2+zy$ in $\P^2$.
Section \ref{sec-example} shows that the
corresponding $U$ has no cohomology but for $H^0(U;\C)\cong \C$.
Hence $H^0(Y;\C)\cong \C$ corresponding to the second
Chern class on $\P^2$, $H^2(Y;\C)\cong\C$ corresponding to the first Chern
class on $\P^2$, and all other
cohomology groups of $Y$ vanish. 
\end{ex}
\begin{ex}
\mylabel{ex-elliptic}
We consider the curve $C=\var(x^2y+y^2z+z^2x)$ in $\P^2$. $\P^2$ is
covered by $U_1=\spec\C[y/x,z/x]$, $U_2=\C[x/y,z/y]$ and
$U_3=\C[x/z,y/z]$.  
We take as coordinates $s=x/z,t=y/z$ on $U_3$,
$U_{1,3}$, $U_{2,3}$ and $U_{1,2,3}$; $s=y/x,t=z/x$ on $U_1$ and
$U_{1,2}$; $s=x/y,t=z/y$ on $U_2$. 
We have the following dehomogenizations for $f$:
\begin{verbatim}
R=QQ[s,t]
f1=s+s^2*t+t^2
f2=s^2+t+s*t^2
f3=s^2*t+t^2+s
\end{verbatim}
Moreover, on two- and threefold intersections $U\cap P_I$ is defined
by the nonvanishing of 
\begin{verbatim}
f12=(s+s^2*t+t^2)*s
f13=(s^2*t+t^2+s)*s
f23=(s^2*t+t^2+s)*t
f123=(s^2*t+t^2+s)*s*t
\end{verbatim}
With \Mtwo\ one computes with the 
command 
\begin{verbatim} 
deRhamAll(g)
\end{verbatim} 
the following cohomology generator table
where $g$ is one of
$f_1,\ldots,f_{1,2,3}$.

\small 
\[
\begin{array}{|cc|c|c|c|c|c|}
\hline
U_1&H^1&2st+1&&&&\\
&&s^2+2t&&&&\\\cline{3-7}
&H^2&1&t&s&&\\\hline
U_2&H^1&\frac{t^2+2s}{2}&&&&\\
&&\frac{2st+1}{2}&&&&\\\cline{3-7}
&H^2&1&t&s&&\\\hline
U_3&H^1&2st+1&&&&\\
&&s^2+2t&&&&\\
&H^2&1&t&s&&\\\hline
U_{1,2}&H^1&2t^2+s&\frac{s^2t+t^2+s}{2}&&&\\
&&-s^3-2st&0&&&\\\cline{3-7}
&H^2&t&t^2&1&s^2&\\\hline
U_{1,3}&H^1&-2t^2-s&s^2t+t^2+s&&&\\
&&s^3+2st&0&&&\\
&H^2&t&t^2&1&s^2&\\\hline
U_{2,3}&H^1&0&\frac{-2st^2-t}{4}&&&\\
&&-s^2t-t^2-s&\frac{-t^2+s}{4}&&&\\\cline{3-7}
&H^2&1&t^2&st&t&\\\hline
U_{1,2,3}&H^1&t^2s^2-t^3&\frac{2t^3+st}{4}&-\frac{s^2t^2+t^3+st}{2}&&\\
&&s^3t+2st^2&\frac{-st^2+s^2}{4}&0&&\\\cline{3-7}
&H^2&t^2&s&t&s^2t&1\\\hline
\end{array}
\]
\normalsize
In this table the generators for $H^1(U_I;\C)$ 
correspond to  columns where the elements of the
top row have to be multiplied with $\frac{ds}{f_I}$ and those of 
the  bottom with
$\frac{dt}{f_I}$. So for example $H^1(U_{1,2,3};\C)$ has three
generators, the first of which is
$\frac{(t^2s^2-t^3)ds+(s^3t+2st^2)dt}{f_{1,2,3}}$. Similarly, the
polynomials listed next to $H^2$ are to be multiplied with
$\frac{ds\,dt}{f_I}$ and then are generators for $H^2(U_I;\C)$. 
So  for example $H^2(U_3;\C)$ has three
generators the last of which is $\frac{s\,ds\,dt}{f_3}$.

We denote these classes by $e_{I:k}$ ($H^1$-generator in column $k$)
and $t_{I:k'}$ ($H^2$-generator in column $k'$) where
$I$ is the index
of the open set in question (for example, $\{1,3\}$ for $U_{1,3}$). 
Thus,
$t_{1,2,3:2}$ is the class $\frac{s\,ds\,dt}{f_{1,2,3}}$. 

As always for a connected set,  the group $H^0(U;\C)$ is a
one-dimensional vector space, and it is here  
generated by the cocycle $(1,1,1)$ of 0-forms. The 0-forms make no
further contribution to the cohomology of $U$.

Using the transformation rules (for example from $U_1$ to
$U_{1,3}$ they say $s\to s^{-1}t, ds\to -s^{-2}t\,ds+s^{-1}\,dt$ and
$t\to s^{-1}, dt\to s^{-2}ds$) one computes that 
\small
\begin{eqnarray*}
e_{3:1}&=&e_{1,3:1}+2e_{1,3:2}=-e_{2,3:1}-4e_{2,3:2},\\
e_{1:1}&=&e_{1,3:1}-e_{1,3:2}=-e_{1,2:1}+4e_{1,2:2},\\
e_{2:1}&=&e_{2,3:1}-2e_{2,3:2}=-\frac{1}{2}e_{1,2:1}-e_{1,2:2},\\
e_{1,3:1}&=&e_{1,2,3:1}+2e_{1,2,3:3},\\
e_{1,3:2}&=&-2e_{1,2,3:3},\\
e_{2,3:1}&=&-e_{1,2,3:1}-4e_{1,2,3:2}-2e_{1,2,3:3},\\
e_{2,3:2}&=&e_{1,2,3:2}+e_{1,2,3:3},\\
e_{1,2:1}&=&e_{1,2,3:1}+8e_{1,2,3:2}+4e_{1,2,3:3},\\
e_{1,2:2}&=&\frac{1}{2}e_{1,2,3:1}+2e_{1,2,3:2}+2e_{1,2,3:3}.
\end{eqnarray*}
\normalsize
This shows that the 1-forms form a complex of vector spaces 
with entries of dimensions
$1+1+1$, $2+2+2$ and $3$ where the matrices have ranks
3 and 3 respectively. Hence the 1-forms make no contribution to the
cohomology of $U$.

Finally, applying the conversion rules to the 2-forms on the various
$U_I$ one obtains a complex that has three entries (on one-, two- and
threefold intersections of open sets) with matrices
$M_{2:1}:\C^{3\times 3}\to \C^{3\times 4}$ and $M_{2:2}:\C^{4\times
3}\to\C^{1\times 5}$. These matrices turn out to have ranks 7 and 
5
respectively. Hence the 2-forms contribute 2 generators of $H^2(U;\C)$
and nothing to $H^3(U;\C)$ or  $H^4(U;\C)$. 

So $U$ has cohomology only in degrees 0 and 2, and $H^2(U;\C)$ is of
dimension 2. It is noteworthy that the first Chern class of $U$ must
be torsion, because it is the pullback from $\P^2$ of a 1-cocycle of
1-forms and we saw that 1-forms make no contribution to the
rational cohomology of $U$. 

From the
Alexander duality exact sequence we see that the complementary curve
$C$ has either Betti numbers 1,2,1 or 1,1,0. The latter case can occur
only if the first Chern class of $U$ is nonzero. Since we know it
vanishes, $C$ has a two dimensional $H^1$ and a one-dimensional
$H^2$. As one can check on a local chart, $C$ is smooth and therefore
topologically $S^1\times S^1$.
\end{ex}

\subsection{Compact Cohomology}
If $Y$ is an affine variety $Y\subseteq X=\A^n$,
then one can compute the cohomology groups with compact support from
Alexander duality on the affine space.
\begin{ex}
Let $f=x^3+y^3+z^3$. Then the de Rham cohomology groups on
$U=\C^3\setminus \var(f)$ have dimensions 1, 1, 2 and 2, which we
computed by \Mtwo\ with
\begin{verbatim}
R=QQ[x,y,z]
deRhamAll(x^3+y^3+z^3)
\end{verbatim}
From Alexander
duality one concludes that (since $\C^3$ has real dimension 6
and is contractible) 
\[
(H^i_c(Y;\C))^*\cong H^{6-i}_Y(X;\C)\cong H^{6-i-1}(U;\C)
\]
for $i<5$ and $H^i_c(Y;\C)=0$ for $i>4$. 
Thus the cohomology groups $H^i_c(Y;\C)$ 
with compact support of $Y$ have dimensions
0, 0, 2, 2, 1 for $i=0,\ldots,4$ and are zero otherwise. 
\end{ex}
One can push the computations a little further in nice situations.
\begin{ex}
If $Y=\var(f)\subset \P^n$ is smooth, then the (usual) 
cohomology of $\var(f)\cap
P_0$ (an affine chart of $\var(f)$) can be computed. For example, if
$f=x^2+yz$, consider the closed subset $Z$ of $Y$ given by $x=0$. This is a
2 point set. Let us now compute the cohomology of $Y\cap
P_0=Y\setminus Z=\var(1+yz)\subset \A^2$. By Alexander duality on $Y$, 
$H^{2\cdot
1-i}(Z;\C)=H^i_Z(Y;\C)^*$. Set $V=\P^2\setminus Z$ and
$U=\P^2\setminus Y$. Then the long exact sequence of sheaf
cohomology on $Y$ gives (with duality incorporated)
\begin{eqnarray*}
&&0\to H^2(Z;\C)^*\to H^0(Y;\C)\to H^0(V\cap Y;\C)\\
&&\phantom{0 }\to H^1(Z;\C)^*\to H^1(Y;\C)\to H^1(V\cap Y;\C)\\
&&\phantom{0 }\to H^0(Z;\C)^*\to H^2(Y;\C)\to H^2(V\cap Y;\C)\to 0.
\end{eqnarray*}
The map $H^{2-k}(Z;\C)^*\to H^k(Y;\C)\cong H^{2-k}(Y;\C)^*$
(by Poincar\'e duality) 
is induced by Alexander duality and really should be
thought of as the dual of the map 
$H^k(Y;\C)\to H^k(Z;\C)$ induced by $Z\into
Y$. 
Alexander duality shows that $H^{2-k}(Z;\C)^*\to H^{2-k}(Y;\C)^*$ is
equivalent to 
$H^{4-2+k}_Z(\P^2;\C)\to H^{4-2+k}_Y(\P^2;\C)$ (i.e., Alexander duality
gives a quasi-isomorphism of the 2-term sequences).  

Consider the commutative diagram
\small
\[
\begin{diagram}
H^{2n-k}(\P^n;\C)&\rTo&H^{2n-k}(V;\C)&\rTo&H^{2n-k+1}_Z(\P^n;\C)&\rTo&H^{2n-k+1}(\P^n;\C)&\rTo&H^{2n-k+1}(V;\C)\\
\dTo_{=}&&\dTo_{\rho^{2n-k}_{U,V}}&&\dTo&&\dTo_{=}&&\dTo_{\rho^{2n-k+1}_{U,V}}\\ 
H^{2n-k}(\P^n;\C)&\rTo&H^{2n-k}(U;\C)&\rTo&H^{2n-k+1}_Y(\P^n;\C)&\rTo&H^{2n-k+1}(\P^n;\C)&\rTo&H^{2n-k+1}(U;\C)
\end{diagram}\]
\normalsize
Inspection shows that $\ker(H^{k+1}_Z(\P^n;\C)\to
H^{k+1}_Y(\P^n;\C))$ is isomorphic to 
\begin{eqnarray*}
\frac{\ker(H^k(V;\C)\to H^k(U;\C))}
{\im(H^k(\P^n;\C)\to H^k(V;\C))\cap \ker(H^k(V;\C)\to H^k(U;\C))}
\end{eqnarray*}
(and zero in the case $k=0$).
This dimension can
be computed by our methods on the level of differential forms, 
by repeatedly applying Algorithm
\ref{alg-lemma-proof}.
It follows that we can evaluate the
dimensions of the kernel and cokernel of $H^{2-k}_Z(\P^2;\C)\to
H^{2-k}_Y(\P^2;\C)$ and hence the dimensions of $H^k(V\cap Y;\C)$. In
our example, $U$ has no nontrivial 
cohomology as pointed out, $V=\P^2\setminus$
two points. We get the following table of dimensions of cohomology
groups:
\[
\begin{array}{|c|c|c|c|c|c|}
\hline
&k=0&k=1&k=2&k=3&k=4\\\hline
H^k(\P^2;\C)&1&0&1&0&1\\\hline
H^k(V;\C)&1&0&1&1&0\\\hline
H^k(U;\C)&1&0&0&0&0\\\hline
H^k_Z(\P^2;\C)&0&0&0&0&2\\\hline
H^k(Z;\C)&2&0&0&0&0\\\hline
H^k_Y(\P^2;\C)&0&0&1&0&1\\\hline
H^k(Y;\C)&1&0&0&0&0\\\hline
\ker(H^k(V;\C)\to H^k(U;\C))&0&0&1&1&0\\\hline
\im(H^k(\P^2;\C)\to H^k(V;\C))&1&0&1&0&0\\\hline
H^k(Y\setminus Z;\C)&1&1&0&0&0\\\hline
\end{array}
\]
\end{ex}
\begin{ex}
We continue Example \ref{ex-elliptic} from the previous subsection. There we
found that $C=\var(x^2y+y^2z+z^2x)$ has Betti numbers 1,2 and 1. We
consider now the open set $V$ in $C$ defined by the nonvanishing of
$z$. On the open set $P_3$ of $\P^2$ this set is the cubic curve
defined by $s^2t+t^2+s$. 
It is easy to see that $C$ meets $z=0$ in 2 points, $Z=\{(0,1,0),
(1,0,0)\}$. 

The long exact sequence for the pair $(C,Z)$ gives 
\begin{eqnarray*}
&&0\to H^2(Z;\C)^*\to H^0(C;\C)\to H^0(V;\C)\to\\
&&\phantom{0 }\to H^1(Z;\C)^*\to H^1(C;\C)\to H^1(V;\C)\to\\
&&\phantom{0 }\to H^0(Z;\C)^*\to H^2(C;\C)\to H^2(V;\C)\to 0\\
\end{eqnarray*}
Of course $H^0(V;\C)\cong \C$ and
$H^2(V;\C)$ is zero because $V$ is topologically 
a non-closed surface. Hence the known data imply that $H^1(V;\C)\cong
\C^3$.
\end{ex}
%
%
%

\section{Toric varieties}
The principles outlined in the previous sections also apply to open and
closed sets within smooth toric varieties. We shall demonstrate this
with an example.
\begin{ex}
Let $X$ be the second Hirzebruch surface $F_2$ (see \cite{Fulton}) 
defined by the complete
fan $\Delta$ in the plane whose rays are the vectors $(1,0)$, $(0,1)$,
$(-1,2)$, $(0,-1)$. We denote the 4 maximal cones by $A,\ldots,D$, the
rays by $AB,\ldots,DA$ and the trivial cone by $ABCD$. We write $O_\sigma$
for the ring of regular functions on the affine variety defined by the
cone $\sigma$. One finds
easily that $O_A=k[(x/z)=s_A,(yz^2/w)=t_A]$, $O_B=k[(z/x)=s_B,(yz^2/w)=t_B]$,
$O_C=k[(z/x)=s_C,(w/yx^2)=t_C]$, $O_D=k[(x/z)=s_D,(w/yz^2)=t_D]$.

Let us first compute the cohomology of $X$. (This is of course well
known from combinatorial methods, see \cite{Fulton}.) On each maximal
cone, a complex quasi-isomorphic to the \v Cech-de Rham complex is simply
given by $(\C\cdot 1)$ concentrated in degree zero. The intersections
of neighboring cones lead to spaces isomorphic to $\C\times \C^*$,
so they have a \v Cech-de Rham complex quasi-isomorphic to the complex
$\left(\C\cdot 1\to\C\cdot\frac{df}{f}\right)$ where $f$ is an appropriately
chosen divisor (corresponding to the ray of 
 intersection). For example, the intersection of the cones $B$ and $C$
leads to the divisor $t_C$
and a corresponding de Rham cohomology generator 
$\frac{d(t_C)}{t_C}$. 

The
intersection of cones $A$ and $C$, and $B$ and $D$, and all higher
intersections  are  2-tori
with \v Cech-de Rham complex quasi-isomorphic to the complex 
$\C\cdot 1\to\C\cdot\frac{ds_A}{s_A}\oplus\C\cdot\frac{dt_A}{t_A}
\to\C\cdot\frac{ds_A\,dt_A}{s_At_A}$.

We combine the $4+6+4+1$ \v Cech-de Rham complexes to a complex 
which computes the cohomology of $X$ in terms of \v Cech cochains of
differential forms.

The de Rham cohomology of $X$ is then generated in
degree 0 by the 0-cocycle $(1_A,1_B,1_C,1_D)$ (which means that on each
2-dimensional cone the chosen function is identically 1). 
The group $H^2(X;\C)$ is of rank two and 
generated by the 1-cocycles of 1-forms
$\alpha=(\frac{ds_A}{s_A},\frac{ds_A}{s_A},0,
0,\frac{ds_C}{s_C},\frac{ds_C}{s_C})$ 
and 
$\beta=(\frac{2ds_A}{s_A},\frac{-2dt_A}{t_A},\frac{-dt_A}{t_A},
\frac{dt_C}{t_C},\frac{dt_C}{t_C},0)$  
where these are the 6
components corresponding to the 6 intersections of the 2-cones, ordered
lexicographically. Finally, $H^4(X;\C)$ can be seen to be generated by
the 2-cocycle of 2-forms
$(\frac{ds_A\,dt_A}{s_At_A},\frac{ds_A\,dt_A}{s_At_A},
\frac{ds_A\,dt_A}{s_At_A},\frac{ds_A\,dt_A}{s_At_A})$
on the triple intersections. All other cohomology groups are zero.
\end{ex}
\begin{ex}
Now we consider the cohomology of the complement of the divisor
$f=w-x^2y+z^2y$ in the surface $X$ of the previous example. 
To that end we compute generators for the de
Rham cohomology for the complement on each affine piece determined by
a cone of $\Delta$. This is done by \Mtwo\ and we use the following notation. 
\begin{center}
\small
\begin{tabular}{|c|c|c|c|}\hline
Cone&Variables &Ring&Divisor of $f$\\\hline
$A$&$s_A=x/z$, $t_A=yz^2/w$&$k[s_A,t_A]$&$f_A=1-s_A^2t_A+t_A$\\
$B$&$s_B=z/x$, $t_B=yx^2/w$&$k[s_B,t_B]$&$f_B=1-t_B+s_B^2t_B$\\
$C$&$s_C=z/x$, $t_C=w/yx^2$&$k[s_C,t_C]$&$f_C=t_C-1+s_C^2$\\
$D$&$s_D=x/z$, $t_D=w/yz^2$&$k[s_D,t_D]$&$f_D=t_D-s_D^2+1$\\
$AB$&$s_A=x/z$, $t_A=yz^2/w$&$k[s_A,t_A,s_A^{-1}]$&$f_{AB}=(1-s_A^2t_A+t_A)s_A$\\
$BC$&$s_C=z/x$, $t_C=w/yx^2$&$k[s_C,t_C,t_C^{-1}]$&$f_{BC}=(t_C-1+s_C^2)t_C$\\
$CD$&$s_C=z/x$, $t_C=w/yx^2$&$k[s_C,t_C,s_C^{-1}]$&$f_{CD}=(t_C-1+s_C^2)s_C$\\
$DA$&$s_A=x/z$, $t_A=yz^2/w$&$k[s_A,s_A,t_A^{-1}]$&$f_{AD}=(1-s_A^2t_A+t_A)t$\\
all others&$s_A=x/z$, $t_A=yz^2/w$&$k[s_A,t_A,s_A^{-1},t_A^{-1}]$&$f_{ABCD}=(1-s_A^2t_A+t_A)st$\\\hline
\end{tabular}
\end{center}
\normalsize

In these local variables, we have the following generators for the
cohomology of the various open sets:
\small
\[
\begin{array}{|c|c|c|c|}\hline
&H^0&H^1&H^2\\\hline
A&1&\frac{2s_At_A\,ds_A+(s_A^2-1)dt_A}{f_A}=A_{1,1}&\frac{ds_A\,dt_A}{f_A}=A_{2:1}\\
&&&\frac{s_A\,ds_A\,dt_A}{f_A}=A_{2:2}\\\hline
B&1&\frac{2s_Bt_B\,ds_B+(s_B^2-1)dt_B}{f_B}=B_{1:1}&\frac{ds_B\,dt_B}{f_B}=B_{2:1}\\
&&&\frac{s_B\,ds_B\,dt_B}{f_B}=B_{2:2}\\\hline
C&1&\frac{2s_C\,ds_C+dt_C}{f_C}=C_{1:1}&\\\hline
D&1&\frac{-2s_D\,ds_D+dt_D}{f_D}=D_{1:1}&\\\hline
AB&1&\frac{-(2t_A+2)ds_A-(s_A^3-s_A)dt_A}{f_{AB}}=AB_{1:1}&\frac{ds_A\,dt_A}{f_{AB}}=AB_{2:1}\\
&&\frac{(s_A^2t_A-t_A-1)ds_A}{f_{AB}}=AB_{1:2}&\frac{s_Ads_A}{f_{AB}}=AB_{2:2}\\
&&&\frac{s_A^2ds_A}{f_{AB}}=AB_{2:3}\\\hline
BC&1&\frac{(s_C^2+t_C-1)dt_C}{f_{BC}}=BC_{1:1}&\frac{ds_C\,dt_C}{f_{BC}}=BC_{2:1}\\
&&\frac{2s_Ct_C\,ds_C+t_C\,dt_C}{f_{BC}}=BC_{1:2}&\frac{s_C\,ds_C\,dt_C}{f_{BC}}=BC_{2:2}\\\hline
CD&1&\frac{(2-2t_C)ds_C+s_C\,dt_C}{f_{CD}}=CD_{1:1}&\frac{ds_C\,dt_C}{f_{CD}}=CD_{2:1}\\
&&\frac{(s_C^2+t_C-1)ds_C}{f_{CD}}=CD_{1:2}&\\\hline
AD&1&\frac{-2s_At_A^2\,ds_A-dt_A}{f_{AD}}=AD_{1:1}&\frac{t_A\,ds_A\,dt_A}{f_{AD}}=AD_{2:1}\\
&&\frac{-2s_At_A^2\,ds_A+(-s^2_At_A+t_A)dt_A}{f_{AD}}=AD_{1:2}&\frac{s_At_A\,ds_A\,dt_A}{f_{AD}}=AD_{2:2}\\\hline
ABCD&1&\frac{(-2t_A^2-2t_A)ds_A-s_A\,dt_A}{f_{ABCD}}=ABCD_{1:1}&\frac{ds_A\,dt_A}{f_{ABCD}}=ABCD_{2:1}\\
&&\frac{2s_A^2t_A^2\,ds_A+s_A\,dt_A}{f_{ABCD}}=ABCD_{1:2}&\frac{t_A\,ds_A\,dt_A}{f_{ABCD}}=ABCD_{2:2}\\
&&\frac{2s_A^2t_A^2\,ds_A+(s_A^3t_A-s_At_A)\,dt_A}{f_{ABCD}}=ABCD_{1:3}&\frac{s_At_A\,ds_A\,dt_A}{f_{ABCD}}=ABCD_{2:3}\\
&&&\frac{s_A^2t_A\,ds_A\,dt_A}{f_{ABCD}}=ABCD_{2:4}\\\hline
\end{array}
\]
\normalsize
The map of 0-cochains to 1-cochains of 1-forms 
is given by the following chart.
\tiny
\begin{eqnarray*}
A_{1,1}&\to& (-AB_{1,1}+2AB_{1,2})+(AC_{1,3})-(AD_{1,2})\\\nonumber
B_{1,1}&\to& (-AB_{1,1}+2AB_{1,2})-(BC_{1,1}+BC_{1,2})-(BD_{1,3})\\\nonumber
C_{1,1}&\to& (-AC_{1,1})-(BC_{1,2})+(CD_{1,1}+2CD_{1,2})\\\nonumber
D_{1,1}&\to& (-AD_{1,1})+(BD_{1,2})-(CD_{1,1}).\nonumber
\end{eqnarray*}
\normalsize
On the other hand, the map from 1-cochains to 2-cochains is as follows.
\tiny
\begin{eqnarray*}
AB_{1,1}&\to&
(ABC_{1,1}+ABC_{1,2}-ABC_{1,3})+(ABD_{1,1}+ABD_{1,2}-ABD_{1,3})\\\nonumber
AB_{1,2}&\to&
(\frac{1}{2}ABC_{1,1}+\frac{1}{2}ABC_{1,2})+(\frac{1}{2}ABD_{1,1}+\frac{1}{2}ABD_{1,2})\\\nonumber
AC_{1,1}&\to& (-ABC_{1,1})+(ACD_{1,1})\\\nonumber
AC_{1,2}&\to& (-ABC_{1,2})+(ACD_{1,2})\\\nonumber
AC_{1,3}&\to& (-ABC_{1,3})+(ACD_{1,3})\\\nonumber
AD_{1,1}&\to& (ABD_{1,2})+(ACD_{1,2})\\\nonumber
AD_{1,2}&\to& (ABD_{1,3})+(ACD_{1,3})\\\nonumber
BC_{1,1}&\to& (ABC_{1,1}+ABC_{1,3})+(BCD_{1,1}+BCD_{1,3})\\\nonumber
BC_{1,2}&\to& (ABC_{1,1})+(BCD_{1,1)}\\\nonumber
BD_{1,1}&\to& (ABD_{1,1})-(BCD_{1,1})\\\nonumber
BD_{1,2}&\to& (ABD_{1,2})-(BCD_{1,2})\\\nonumber
BD_{1,3}&\to& (ABD_{1,3)})-(BCD_{1,3})\\\nonumber
CD_{1,1}&\to& (-ACD_{1,2})-(BCD_{1,2})\\\nonumber
CD_{1,2}&\to&
(\frac{1}{2}ACD_{1,1}+\frac{1}{2}ACD_{1,2})+(\frac{1}{2}BCD_{1,1}+\frac{1}{2}BCD_{1,2})\nonumber
\end{eqnarray*}
\normalsize
Finally, the map to 3-cochains is given by the following table.
\tiny
\[
\begin{array}{|c|c|c|c|}
\hline
&ABCD_{1,1}&ABCD_{1,2}&ABCD_{1,3}\\\hline
ABC_{1,1}&     1&0&0\\
ABC_{1,2}&     0&1&0\\
ABC_{1,3}&     0&0&1\\
ABD_{1,1}&     -1&0&0\\
ABD_{1,2}&     0&-1&0\\
ABD_{1,3}&     0&0&-1\\
ACD_{1,1}&     1&0&0\\
ACD_{1,2}&     0&1&0\\
ACD_{1,3}&     0&0&1\\
BCD_{1,1}&     -1&0&0\\
BCD_{1,2}&     0&-1&0\\
BCD_{1,3}&     0&0&-1\\
\hline
\end{array}
\]
\normalsize
Taking cohomology (in \Mtwo\ for example with the commands {\tt ker,
image, gens} and {\tt \%}) we see that $H^2(U;\C)$
is generated by the 1-cocycle of 1-forms
$\alpha=(-AB_{1,2},-\frac{1}{2}AC_{1,1}-\frac{1}{2}AC_{1,2},0,0,\frac{1}{2}BD_{1,1}+\frac{1}{2}BD_{1,2},CD_{1,2})$
(this is the same $\alpha$ that was a cohomology generator on
$X$, now expressed in terms of the cochains on $U$). 
This is the only contribution 1-forms make to the cohomology of $U$.

On the 2-forms we have the following maps for \v Cech cocycles from
single to double intersections:
\tiny
\begin{eqnarray}
\\\nonumber
A_{2,1}\to (AB_{2,2})+(AC_{2,3})+(AD_{2,1)}\\\nonumber
A_{2,2}\to (AB_{2,3})+(AC_{2,4})+(AD_{2,2})\\\nonumber
B_{2,1}\to (AB_{2,2})-(BC_{2,1})+(BD_{2,3})\\\nonumber
B_{2,2}\to (AB_{2,1})-(BC_{2,2})+(BD_{2,2})\nonumber
\end{eqnarray}
\normalsize
The matrix for 2-forms from double to triple intersections can be obtained (we
omit the obvious maps from $AC$ and $BD$ to any triple intersection)
from the equations
\tiny
\begin{eqnarray*}
&AB_{2,1}=BC_{2,2}=-ABCD_{2,2},&\\
&AB_{2,2}=AD_{2,1}=BC_{2,1}=ABCD_{2,3},&\\
&AB_{2,3}=AD_{2,2}=ABCD_{2,4},&\\
&CD_{2,3}=ABCD_{2,3}.&
\end{eqnarray*}
\normalsize
Linear algebra shows that the 2-forms do not contribute to the cohomology of $U$.
Thus, $H^0(U;\C)\cong H^2(U;\C)\cong \C$ and all other cohomology
groups are zero.

It follows from the long exact sequence in \ref{subsec-closed} 
relating the open set and the
variety of $f$ that $\var(f)$ has its cohomology concentrated in
degree 0 and 2 and both are one-dimensional. This is because
$H^2(X;\C)\to H^2(U;\C)$ is not the zero map since $\alpha\not =0$ in
$H^2(U;\C)$. On the other hand,
the cocycle $2\alpha-\beta$ that generates cohomology on $X$ is zero on
$U$: as
\begin{eqnarray*}
\alpha&=&(-AB_{1,2},-\frac{1}{2}AC_{1,1}-\frac{1}{2}AC_{1,2},0,0,\frac{1}{2}BD_{1,1}+\frac{1}{2}BD_{1,2},CD_{1,2}),\\
\beta&=&(-2AB_{1,2},AC_{1,3}-AC_{1,2},AD_{1,1}-AD_{1,2},BC_{1,1},BD_{1,1}+BD_{1,3},0).
\end{eqnarray*}
one sees that 
$2\alpha-\beta=d(-A_{1,1}+B_{1,1}+C_{1,1}+D_{1,1})$. 
Hence $2\alpha-\beta$ is the zero class in $H^2(U;\C)$.
\end{ex}
\bibliography{bib}
\bibliographystyle{abbrv}

\end{document}